\documentclass[12 pt]{article}
\usepackage[final]{epsfig}      
\usepackage{verbatim}
\usepackage{amsfonts,amssymb,latexsym,amscd,theorem,epsfig,graphics}
\newtheorem{theorem}{Theorem}[section]
\newtheorem{lemma}{Lemma}[section]
\newtheorem{remark}{Remark}[section]

\newtheorem{proposition}{Proposition}[section]

\newtheorem{claim}{Claim}[section]
\newtheorem{definition}{Definition}[section]

\begin{document}
\newcommand{\eps}{{\varepsilon}}
\newcommand{\proofend}{$\Box$\bigskip}
\newcommand{\C}{{\mathbf C}}
\newcommand{\Q}{{\mathbf Q}}
\newcommand{\R}{{\mathbf R}}
\newcommand{\Z}{{\mathbf Z}}
\newcommand{\RP}{{\mathbf {RP}}}
\newcommand{\CP}{{\mathbf {CP}}}
\newcommand{\al}{\alpha}

\title {Embeddings With Multiple Regularity}
\author{Gordana Stojanovic\\
{\it Brown University and Penn State University}\\
}
\date{}

\maketitle

\begin{abstract}
We introduce $(k,l)$-regular maps, which generalize two previously studied classes of maps: affinely $k$-regular maps and totally skew embeddings. We exhibit some explicit examples and obtain bounds on the least dimension of a Euclidean space into which a manifold can be embedded by a $(k,l)$-regular map. The problem can be regarded as an extension of embedding theory to embeddings with certain non-degeneracy conditions imposed, and is related to approximation theory. 
\end{abstract}

\section{Introduction}\

Two lines in a Euclidean space are called \emph{skew} if they are not parallel and do not intersect. 
A submanifold $M^n$ of $\R^N$ is said to be \emph{totally skew} if arbitrary tangent lines to $M$ at any two distinct points are skew. Equivalently, one can define an immersion $f:M^n\to\R^N$ to be \emph {totally skew} if for all $x, y \in M^n$ the tangent spaces $df(T_x M)$ and $df(T_y M)$ as affine subspaces of $\R^N$ have an affine span of maximal possible dimension, that of $2n+1$. 

Totally skew embeddings have been introduced and studied in \cite{G-T}. Other closely related classes of embeddings into affine and projective spaces defined in terms of mutual positions of tangent spaces at distinct points are skew embeddings and T-embeddings and they have also received a considerable amount of attention, see \cite {Gh1,Gh2,G-S,G-T,S-S,S-T,Ta,T-T}. 

Another, seemingly less closely related, class of embeddings are so-called $k$-regular maps, and their affine version, introduced by Borsuk in \cite{Bo}. A continuous map $f:X\to\R^N$ is called \emph{k-regular}  (respectively \emph{affinely $k-1$-regular}) if the images under $f$ of $k$ distinct points are linearly (respectively affinely) independent.\footnote{The definitions of both types of $k$-regularity that we adopt here are those of Handel, see for example \cite{H-S}. The definitions used by other authors are in essence the same, but the role of $k$ differs from one author to another. We adapt all the quoted results to our definition.} 
The study of $k$-regular maps was motivated by the theory of Chebyshev approximation. It was conducted by non-algebro-topological methods in \cite{Bo,BRS}, while Handel \cite{C-H,H1, H2, H3, H4, H-S} introduced cohomological methods using configuration spaces to obtain various existence and non-existence results. Vassiliev \cite{V} independently studied $k$-regular maps under the name `$k$-interpolating spaces of functions', using topological methods similar to those of Handel. He was interested in the interpolating properties of a finite dimensional space of continuous functions on a topological space. Namely, he calls a finite dimensional space $L$ of continuous functions on a topological space $M$, \emph{$k$-interpolating} if every real-valued function on $M$ can be interpolated at arbitrary $k$ points of $M$ with an appropriate function from $L$. The connection with $k$-regular maps is as follows: the functions $f_1,\ldots,f_N$ span a $k$-interpolating space of functions if and only if the map $f=(f_1,\ldots,f_N)$ is $k$-regular. In other words, $f$ is $k$-regular if and only if we can prescribe values at any distinct $k$ points of $M$ for functions in the span of coordinate functions of $f$.  

One of the main questions that arises in the study of all mentioned maps is to find the lowest possible dimension of the target Euclidean space which allows them. For example, for a given manifold $M^n$, what is the smallest dimension $N=N(M^n)$ such that $M^n$ admits a totally skew embedding in $\R^N$? As is, this question has been answered for very few manifolds. Results are available only for line, circle and plane: $N(\R^1)=3, N(S^1)=4$, $N(\R^2)=6$, see \cite{G-T}. Ghomi and Tabachnikov actually give totally skew embeddings of line, circle and plane in the Euclidean space of minimal possible dimension and these are the only known explicit examples of optimal totally skew embeddings. According to the same authors \cite{G-T}, dimension $n$ submanifolds of $\R^N$ are generically totally skew when $N\geq 4n+1$. This abundance of totally skew embeddings contrasted with the scarcity of available examples points to another object of investigation: finding more of them.

The same question can be asked for $k$-regular maps, and one result, that both Handel \cite{H4} and Vassiliev \cite{V} reached, is for instance, that when $k$ is even, $N(S^1)=k+1$, and when $k$ is odd, $N(S^1)=k$. While the result for odd $k$ is almost immediate, to achieve the result for $k$ even, they both used nonelementary topological methods, in particular, characteristic classes.

We introduce a class of regular maps, so called $(k,l)$-regular maps, which generalize both totally skew embeddings and affinely $k$-regular maps, and ask the same question of determining minimal dimensional target Euclidean space. This problem can be regarded as an extension of investigations that led to the birth of embedding theory - an extension to the embeddings with certain prescribed non-degeneracy conditions. The interpretation of $(k,l)$-regular maps in the language of the approximation theory is as follows: it turns out that a smooth map $f=(f_1,\ldots,f_N):M^n\to \R^N$ on a smooth manifold $M^n$ is $(k,l)$-regular if and only if for every function in the span of $1,f_1,\ldots,f_N$ we can prescribe not only values at any distinct points $x_1,\ldots,x_k,y_1,\ldots,y_l$ but directional derivatives as well in any direction at the last $l$ points. Thus, the existence of $(k,l)$-regular maps is equivalent to the possibility of interpolating functions on $M^n$ through any $k+l$ points and up to the first order derivatives in arbitrary directions at the last $l$ points. Finally, let us mention that there is an obvious connection with recent work of Arnold and his school, see \cite{A}.
 
In this paper, we generalize existing estimates for totally skew and affinely $k$-regular maps to our class, provide explicit examples in the case of line, circle and plane and determine the minimal target spaces for curves. We only employ non-algebro-topological methods, which leaves plenty of room for further investigations in the topology of these embeddings. 
  
\smallskip

{\bf Acknowledgments}. I would like to thank my advisor Sergei Tabachnikov for suggesting this problem as well as for his constant support and guidance. I would also like to thank Mohammad Ghomi and Bruce Solomon for interesting discussions and Anatole Katok for continual support during my stay at Penn State University, where this work has been carried out. 

\section {Definition of \boldmath$(k,l)$-regular maps}
%\begin{definition} Affine subspaces $V_1,\ldots,V_k\subset\R^N$ of dimensions $n_1,\ldots, n_k$ respectively, are said to be in general position if their affine span has the dimension $$(n_1+1)+\ldots+(n_k+1)-1.$$ 
%\end{definition}

We will start with the definition of affine independence. There are many different (albeit equivalent) ways to define affine independence, but we settle with the following one.

\begin{definition} Affine subspaces $V_1,\ldots,V_k\subset\R^N$ are said to be \emph{affinely independent} if their affine span has the maximal possible dimension that the affine span of affine spaces of respective dimensions can have in any given affine ambient space.
\end{definition} 

For example, the affine span of two lines may have  dimension 1, 2 or 3 depending on their position. The maximal of these is three dimensional, and so any two lines are affinely independent if their affine span is three dimensional. Thus no two lines in $\R^2$ are affinely independent. %, even though their span may be all of $\R^2$ (which is maximal possible within the ambient space). 

From now on, when we say that some span is maximal possible we will mean maximal possible regardless of the ambient space. In this terminology, no two lines in $\R^2$ will have maximal possible affine span.

One can calculate that if $V_1,\ldots,V_k\subset\R^N$ have dimensions $n_1,\ldots, n_k$ respectively than they are affinely independent if and only if their affine span has the dimension $$(n_1+1)+\ldots+(n_k+1)-1.$$  
 
\begin {definition} Let $M^n$ be an $n$-dimensional manifold. Let $k$ and $l$ be non-negative integers, not both equal to 0. We will call a smooth map \mbox{$f:M^n\to\R^N$} \emph{$(k,l)$-regular} if for every set of distinct  points $x_1,\ldots,x_k,y_1,\ldots,y_l$ of $M$ and of $l$ tangent lines $L_i\subset T_{y_{i}} M, i=1,\ldots,l$ the set of points and lines
$$ f(x_1),\ldots,f(x_k),df(L_1),\ldots,df(L_l)$$ is affinely independent. 
\end{definition}

\begin{remark}{\rm When $l=0$, the notion of $(k,l)$-regularity coincides with the notion of affine $k-1$-regularity as defined in \cite {H-S}. On the other hand, a smooth map $f$ is (0,2)-regular if and only if it is totally skew. Thus, the notion of \mbox{$(k,l)$-regularity} generalizes both totally skew and (affinely) $k$-regular maps.}
\end{remark}

We will often say that a manifold is $(k,l)$-regular with the understanding that it is a submanifold of some Euclidean space, and the inclusion map is \mbox{$(k,l)$-regular}.

\begin{definition}
The $(0,l)$-regular maps we will also call $l$-totally-skew.
\end{definition}

\begin{remark}\label{obvious}{\rm The following properties are more or less obvious:
\begin{itemize}
\item If a map is $(k,l)$-regular then it is also $(k',l')$-regular, for all $k'\leq k$, $l'\leq l$.
\item Any (2,0)-regular map is one-to-one. 
\item Any 1-totally-skew map is an immersion.
\item Any 2-totally-skew map is nothing else but a totally skew embedding.
\item Every restriction of a $(k,l)$-regular map is $(k,l)$-regular.
\end{itemize}}
\end{remark}
\smallskip

One may be tempted to adopt a more general definition 
of regularity: instead of requiring affine independence of a collection of points and lines, one may consider arbitrary dimensional vector subspaces of tangent spaces at certain number of distinct points of $M$ and require the affine independence of their images under $df$. However, it turns out that this doesn't add any generality to the definition, as stated in the following lemma, which we leave without proof. 

\begin{lemma} \label{reduction}
Let $M^n$ be an $n$-dimensional manifold and fix some positive integers $n_1,\ldots,n_l$ not greater than n. Then, the map $f:M^n\to\R^N$ is $(k,l)$-regular if and only if for arbitrary distinct points \mbox{$x_1,\ldots,x_k,y_1,\ldots,y_l\in M$} and arbitrary vector spaces $V_i\subset T_{y_{i}}M,~ dim(V_i)=n_i,~i=1,\ldots,l$ the affine spaces $$f(x_1),\ldots,f(x_k),df(V_1),\ldots,df(V_l)$$ are affinely independent. 
\end{lemma}

%We note here that not only do $k,l$-regular maps generalize affinely $k$-regular maps and $l$-totally skew maps, but in some sense they interpolate between them. Namely, $k$-regular maps and $l$-totally skew maps are on the opposite ends of the spectrum of $k,l$-regular maps. On the one hand, for $k$-regular maps we consider the affine span of $k$ points of $f(M)$, that is, of $k$ 0-dimensional subspaces of tangent spaces to $f(M)$. On the other hand, for $l$-totally skew maps we consider the affine span of $l$ lines tangent to $f(M)$, or,which is equivalent by previous lemma, the affine span of full tangent spaces to $f(M)$. \

To fix notation we will denote $N_{k,l}=N_{k,l}(M^n)$ to be the least possible integer $N$, such that $M^n$ admits a $(k,l)$-regular map into $\R^N$.
 
\begin {remark} {\rm A simple dimension counting argument shows that any $n$-dimensional $(k,l)$-regular manifold requires at least $k+(n+1)l-1$ dimensions of ambient space, immediately yielding that 
\begin{equation} \label{count}
N_{k,l}(M^n)\geq k+(n+1)l-1.
\end{equation} 
In particular, 
\begin{equation} \label{count1}
N_{k,l}(\R)\geq k+2l-1.
\end{equation}}
\end{remark}
\begin{remark}{\rm 
Since manifolds are locally Euclidean, it follows that 
$$ N_{k,l}(M^n)\geq N_{k,l}(\R^n). $$}
\end{remark}

The following observation will be useful in the proofs. 

\begin{lemma} \label {notreg}
{\rm An embedding $f:M^n\to\R^N$ is \emph{not} $k,l$-regular if and only if there exists an \mbox{$(k+(n+1)l-2)$-dimensional} affine subspace that touches $M$ at $l$ points and intersects it in additional $k$ points.} 
\end{lemma}

\noindent {\it Proof.} It follows from Lemma \ref{reduction} where $V_j$'s are taken to be full tangent spaces.
\proofend.

Henceforth we shall assume that $f$ is always an embedding, simply because in most interesting cases  (e.g. whenever $k+l>1$), $(k,l)$-regular maps actually are embeddings.

\section{Examples and bounds}

The following proposition, together with the last item in Remark \ref{obvious} and the existence of $(k,l)$-regular embeddings of real line (Thm \ref{curves}) implies the existence of $(k,l)$-regular embeddings of any manifold into some Euclidean space, and as a consequence, $N_{k,l}$ is well defined.

\begin{proposition} \label {tensor}
Let $M$ and $N$ be $(k,l)$-regular submanifolds in  euclidean spaces $U$ and $V$ respectively. Then, the  embedding $$f:M \times N \to (U\otimes V) \oplus (U\otimes \R)\oplus (\R \otimes V),~ f(x,y)=(x \otimes y, x\otimes1,1\otimes y)$$ is $(k,l)$-regular. 
\end{proposition}

The next theorem improves on the estimate (\ref{count}) and provides an upper bound as well. The lower bound is obtained by extending the argument in \cite{BRS}, while the upper bound is obtained similarly as in the case of totally skew embeddings \cite{G-T}.\

\begin{theorem} \label{main}
For any manifold $M^n$,
$$\Biggl[\frac{k}{2}\Biggr]n+\biggl[\frac{k-1}{2}\biggr]+(n+1)l \leq N_{k,l}(M^n)\leq (n+1)k+(2n+1)l-1.$$Moreover, generically any submanifold $M^n$ in $R^{(n+1)k+(2n+1)l-1}$ is $(k,l)$-regular.
\end{theorem}

Setting $k=0$ and $l=2$ gives us the theorem for totally skew embeddings obtained in \cite{G-T}. On the other hand, the upper bound for the opposite case, $k=2$ and $l=0$, says simply that $M^n$ always embeds in $R^{2n+1}$, a well known fact. The theorem also generalizes the lower bound for affinely $k$-regular maps given in \cite {BRS}.

When the manifold is closed, we have a better (by 1) lower bound in the case of $l$-totally-skew embeddings.

\begin {theorem} \label{closed}
Let $M^n$ be a closed manifold. Then, $$N_{0,l}(M^n)\geq (n+1)l.$$
\end{theorem}

Finally, we give examples of $(k,l)$-regular embeddings of real line, circle and plane. In the case of curves, these embeddings happen to be optimal, deciding $N_{k,l}$ when $n=1$. 

\begin{proposition} \label{examples}
\begin{enumerate}
\item The map $\gamma:\R \to \R^{k+2l-1}$ given by $$\gamma(t)=(t,t^2,\ldots,t^{k+2l-1})$$ is $(k,l)$-regular.
\item The map $\gamma:S^1\to\R^{2k'+2l}$, given by $$\gamma(\al)=( \cos \al, \sin \al,\cos 2\al,\sin 2\al,\ldots,\cos (k'+l)\al,\sin (k'+l)\al)$$ is $(k,l)$-regular for $k=2k'+1$ (and hence for $k=2k'$ as well).
\item The map $\gamma:\R^2 \cong \C \to \C^{k+2l-1} \cong \R ^{2(k+2l-1)}$ given by $$\gamma(z)=(z,z^2,\ldots,z^{k+2l-1})$$ is $(k,l)$-regular.
\end {enumerate}
\end{proposition}

\begin{theorem} \label {curves} One has:
\begin{enumerate}
\item $N_{k,l}(\R)= k+2l-1$.\
\item $N_{k,l}(S^1)= \left \{ \begin{array}{ll}
k+2l, & \textrm {$k$ is even}\\
k+2l-1, & \textrm{$k$ is odd}.
\end{array} \right.$
\end{enumerate}
\end{theorem}

\section{Proofs}\
Let $i$ be the embedding of $\R^N$ into $\R^{N+1}$ as the height 1 hyperplane, that is let $$i:\R^N\to \R^{N+1},~i:x\mapsto (x,1).$$ Let $AG_n(N)$ denote the (affine Grassmanian) manifold of $n$-dimensional affine subspaces of $\R^N$, and $G_{n+1}(N+1)$ be the (Grassmanian) manifold of the ($n+1$)-dimensional subspaces of $\R^{N+1}$. Then, $i$ induces the canonical embedding $$AG_n(N)\to G_{n+1}(N+1),$$ given by assigning to each point $p\in\R^N$ the line $\ell (p)\subset\R^{n+1}$ which passes through the origin and $(p,1)$. We will also call this embedding $i$ without any fear of confusion.\

The following fact, which is immediate consequence of the definition, we state as a lemma. 

\begin {lemma} \label{linear}
The affine span of the affine subspaces $V_1,\ldots,V_k\subset\R^N$ is maximal possible if and only if the linear span of the corresponding vector subspaces $i(V_1),\ldots,i(V_k)\subset\R^{N+1}$ maximal possible. The dimension of this linear span is one larger than the dimension of the affine span of $V_1,\ldots,V_k$. 
\end{lemma}
 
\begin{remark} \label {operative}
{\rm We will use notation $\tilde x=(x,1)$ for points $x\in M$, but \mbox{$\tilde u=(u,0)$} for vectors $u$ tangent to $M$. As a consequence of the previous lemma, it follows that $f:M\to \R^N$ is $(k,l)$-regular if and only if for distinct points $x_1,\ldots, x_k,y_1,\ldots, y_l\in M$ and non-zero  vectors $u_j \in T_{y_j}M,j=1,\ldots,l$, the vectors $\tilde{x_1},\ldots,\tilde{x_k},\tilde{y_1},\ldots,\tilde{y_l},\tilde{u_1},\ldots,\tilde{u_l}\in \R^{N+1}$ are linearly independent.}
\end{remark}
\smallskip

\noindent  {\bf Proof of Proposition \ref{tensor}.}\
First note that a tangent vector to $f(x,y)$ is of the form $(u\otimes y+x\otimes v,u \otimes 1, 1 \otimes v)$ where $u\in T_x M$ and $v\in T_y N$. Fix arbitrary distinct points $(x_i,y_i)\in M\times N, i =1,\ldots,k+l$ and arbitrary non-zero vectors 
$(u_j,v_j)\in T_{(x_j,y_j)}M\times N, j=k+1,\ldots,k+l.$ 
%Let $\mathcal M$ and $\mathcal N$ be the set of all distinct $x_i$'s and respectively $y_i$'s that appear in (\ref{points}), and let $\mathcal U$  and $\mathcal V$ be the set of all non-zero $u_j$'s and $v_j$'s respectively that appear in the list (\ref {vectors}). 

 Because of remark \ref{operative}, it suffices to show, after the identification $$(U \otimes V)\oplus(U\otimes\R)\oplus(V\otimes \R)\oplus \R \cong (U\oplus \R)\otimes (V \oplus \R)$$ that the vectors in the set 
$$\mathcal  S=\{\tilde{x_i}\otimes \tilde{y_i}| i=1,\ldots,k+l\} \cup \{\tilde u_j \otimes \tilde y_j + \tilde x_j \otimes \tilde v_j| j=k+1,\ldots,k+l\}$$ 
are linearly independent. 

Since $M$ is $(k,l)$-regular, we know that the vectors in $\mathcal M=\{\tilde {x_i}, \tilde{u_j}|u_j\neq0\}$ are linearly independent and similarly for $\mathcal N=\{\tilde {y_i}, \tilde {v_j}|v_j\neq0\}.$ But now we observe that every nontrivial linear combination of vectors in $\mathcal S$ would also be a nontrivial linear combination of vectors in $$\{e\otimes f| e\in \mathcal M,f \in \mathcal N\},$$ and that is not possible.
\proofend

\noindent  {\bf Proof of Theorem \ref{main}.} First we prove the lower bound.\
We will make use of the following theorem, affine version of which is proved in \cite{BRS}. For one line proof of the even $k$ case of the theorem stated here, see \cite {C-H}.

\begin{theorem}{\rm \bf(Boltyanski-Ryzhkov-Shashkin)}\label{thm}
If a $k$-regular map of $\R^n$ into $R^N$ exists, then $N \geq [\frac{k}{2}]n+[\frac{k+1}{2}]$.
\end{theorem} 

Since  the inequality holds trivially when $k+l=1$, and $f$ is an embedding when $k+l>1$, we can assume without loss of generality that $M^n$ is a \mbox{$(k,l)$-regular} submanifold of $R^N$.

%For $x \in \R^N$ denote by $\tilde x= (x,1)$. For $v \in T_y M\subset T_y \R^N\cong \R^N$ , we denote $\tilde v=(v,0).$ Let $x_1,\ldots,x_k,y_1,\ldots,y_l$ be arbitrary distinct points in $M$ and $v_j \in T_{y_j} M$. Then, by lemma \ref{linear}, $f$ is $(k,l)$-regular if and only if the vectors $\{\tilde x_i\},\{ \tilde y_j\},\{\tilde v_j\}$  with $i=1,\ldots,k, j=1,\ldots,l$ are linearly independent and this is equivalent to 
%\begin{displaymath}
%|a|^2+|b|^2+|v|^2\neq 0 \Rightarrow \sum_{i=1}^k a_i\tilde x_i + \sum_{j=1}^l b_j\tilde y_j + \sum_{j=1}^l \tilde v_j \neq 0,
%\end{displaymath}
%where $a=(a_1,\ldots,a_k) \in \R^k,~ b=(b_1,\ldots,b_l) \in \R^l$ are arbitrary, $v=(v_1,\ldots,v_l) \in \R^{nl}$ with $v_j \in T_{y_j}M \cong \R^n$  and $|\cdot|$ is the standard Euclidean norm.\
 
Without loss of generality, assume $l>0$; otherwise, there is nothing to prove. Now fix some arbitrary distinct points $y_1,\ldots,y_l\in M^n$, and let $W$ be the affine span of $T_{y_1}M ,\ldots, T_{y_l}M$. The inclusion $M\to \R^N$ induces a map $M^n- \{y_1,\ldots,y_l\} \to \R^{N+1} /i(W)\cong \R^{N+1-(n+1)l}$ . Since $M^n$ is $(k,l)$-regular, the induced map is $k$-regular. Because of Theorem \ref{thm} we have that $$N+1-(n+1)l \geq \Biggl[\frac{k}{2}\Biggr]n + \Biggl[\frac{k+1}{2}\Biggr],$$ which proves the lower bound. 

%\begin{displaymath}
%W=\bigg\{\sum_{j=1}^l b_j\tilde y_j + \sum_{j=1}^l \tilde v_j~|~b_j\in\R, v_j\in T_{y_j}M\cong\R^n\bigg\} 
%\end{displaymath}
\
To prove the upper bound, the strategy is to first embed the manifold $M$ in a large Euclidean space $R^N$ as a $(k,l)$-regular submanifold. This can be done, for instance, using Proposition \ref{tensor}. Then, one reduces the dimension of the target space by succesive projections to hyperplanes all the while preserving $(k,l)$-regularity.

To do that, we  project to any hyperplane, centrally from a point $A$ outside of it. The projected manifold will still be $(k,l)$-regular if we choose $A$ away from $U$, the union of affine spans of all $k+l$-tuples of points of $M$ and $l$ tangent spaces at the last $l$ points. Since $U$ consists of points
\begin{displaymath} 
\sum_{i=1}^k a_i x_i+ \sum_{j=1}^l b_j y_j+\sum_{j=1}^l v_j
\end{displaymath} 
where $$a_i,b_j\in\R,x_i,y_j\in M^n, v_j\in T_{y_j}M^n,
\sum_{i=1}^k a_i+\sum_{j=1}^l b_j=1$$
one can calculate using local charts that dim $U=k(n+1)+l(2n+1)-1$.  
Thus, for $N>k(n+1)+l(2n+1)-1$ one can find $A$ not in $U$ and reduce $N$ by 1. This proves the upper bound. 

To prove that the upper bound is generically true, that is, that every embedding $M^n \to \R^{k(n+1)+l(2n+1)-1}$ can become $(k,l)$-regular after an arbitrarily small perturbation, one uses Thom's transversality theorem. The proof is exactly the same as in \cite {G-T}, so we omit it.  
\proofend

\smallskip

\noindent  {\bf Proof of Theorem \ref{closed}.}\
Thm \ref{closed} follows immediately from the next proposition and remark \ref{notreg}.
\begin{proposition}
If $M^n$ is a closed submanifold of $\R^{(n+1)l-1}$, then there exists a hyperplane tangent to $M$ at $l$ distinct points.
\end{proposition}

\noindent {\it{Proof.}}\
%In this proof we will use standard facts from convex geometry, which can be found in any introductory book on the subject, see for example (quote Barvinok).\ 
Without loss of generality, we may assume that $M=M^n$ does not belong to an affine hyperplane in $\R^{(n+1)l-1}$; if it does, we are done. Let us denote by $\mathrm{C}(M)$ the convex hull of $M$. 

\begin{claim} \label{convex} There exists a point $x\in\partial\mathrm{C}(M)$ that is not a convex combination of $l-1$ or fewer points of $M$.
\end{claim}

\textit{Proof.} 
The set of all convex combinations of all $l-1$-tuples of points of $M$,
$$ S = \{x= a_1 y_1+\ldots+a_{l-1} y_{l-1}|\sum_{i=1}^{l-1} a_i=1,~ a_i\geq 0, y_i\in\ M, i=1,\ldots,l-1\},$$ 
is a set of dimension $(n+1)l-2-n$, containing $M$. On the other hand, since $M$ does not belong to any hyperplane, neither does $\mathrm{C}(M)$. Therefore, $\mathrm{C}(M)$ is a convex set of dimension $(n+1)l-1$, its boundary $\partial\mathrm{C}(M)$ has dimension $(n+1)l-2$. Thus, there must exist a point $x$ in $\partial\mathrm{C}(M)$, but not in $S$, as claimed.   

\smallskip

Consider now the support hyperplane $H$ at $x$. It is a hyperplane through $x$, so that $\mathrm{C}(M)$ is contained in a closed half-space bounded by $H$. Since $M$ is closed, $\mathrm{C}(M)$ is compact, and therefore $x \in \mathrm{C}(M)$. By the theorem of Carath\'eodory, see for example \cite{Ba}, every point in $\mathrm{C}(M)$ is a convex combination of \emph{at most} $(n+1)l$ points, so 
$$x=\ a_1 y_1+\ldots+a_{(n+1)l} y_{(n+1)l},~~~\mathrm{where}~~~ \sum_{i=1}^{(n+1)l} a_i=1,~ a_i\geq 0$$ and $y_i\in\ M$ for $i=1,\ldots,(n+1)l$. \ 
%modify or perhaps omit the following proof.
%\begin{claim} If $a_i\neq 0$, then $y_i\in H$.
%\end{claim}
%\textit{Proof.} All the points $y_i$ can be grouped in two sets: those in $H$, and those in an \emph{open} half-space bounded by $H$. By reordering our points if necessary we may write
%$$ x=\sum_{i=1}^m a_i y_i + \sum_{i=m+1}^{2k}a_i y_i,$$where $y_i\in H, i\leq m$ and $y_i\notin H, i>m$ 
%for some m, $1\leq m\leq {2k}$. Setting $t=\sum_{i=1}^m a_i,$ we have $1-t=\sum_{i=m+1}^{2k}a_i$, and we can rewrite this as $$x=t z+ (1-t)w, ~\mathrm {where} ~ z\in H~\mathrm{and}~ w\notin H~\mathrm{and}~ 0\leq t\leq 1.$$ If $t$ were not 1, $x$ would not be in $H$. Thus, $
%1-t=\sum_{i=m+1}^{2k}a_i=0$ and $a_i\geq 0$ imply that $$x=\sum_{i=1}^m a_i y_i, ~\mathrm{with}~y_i\in H.$$\
However, all the points $y_i$ with non-zero coefficients $a_i$ in the above convex combination must belong to $H$, because otherwise $x$ wouldn't be in $H$. So, it is exactly at those points where $H$, our support hyperplane, touches $M$. And by claim \ref{convex} there must be at least $l$ of them. Which is exactly what we wanted to prove.
\proofend

\noindent {\bf{Proof of Proposition \ref{examples}}}\

\noindent {\it{1.The open case.}}\
\smallskip

We will argue by contradiction. Suppose that the curve 
 $$\gamma(t)=(t,t^2,\ldots,t^{k+2l-1}),~ t\in\R$$ is not $(k,l)$-regular. Then there exists a hyperplane $H$,
\begin{displaymath}
a_0+\sum_{i=1}^{k+2l-1}a_i x_i=0,
\end{displaymath} 
tangent to the curve $\gamma(t)$ at $l$ distinct points and intersecting it in at least $k$ additional points. But this means that the polynomial 
$$f(t)=a_0+\sum_{i=1}^{k+2l-1}a_i t^i,$$ of degree $k+2l-1$, has $l$ double and $k$ simple roots. Contradiction.
\smallskip

\noindent{\it{2.The closed case.}}\
\smallskip

Again, suppose that the curve $\gamma:S^1\to\R^{2k'+2l}$, $$\gamma(\al)=(\cos \al, \sin \al, \cos 2\al, \sin 2\al,\ldots, \cos (k'+l)\al, \sin (k'+l)\al)$$ is not $(k,l)$-regular for $k=2k'+1$. Then, just as in the open case we obtain %a hyperplane $H$,
%\begin{displaymath}
%a_0+\sum_{i=1}^{k'+l}a_i x_i+\sum_{i=1}^{k'+l}b_i %y_i=0,  
%\end{displaymath}
%having $l$ tangencies and $k$ aditional intersections with $\gamma$, and) 
a function
\begin{displaymath}
f(\al)= a_0+\sum_{i=1}^{k'+l}a_i \cos i\al +\sum_{i=1}^{k'+l}b_i \sin i\al.
\end{displaymath}
having $2k'+2l+1$ roots on the interval $[0,2\pi)$ when counted with multiplicites. However, $f(\al)$ is a trigonometric polynomial of degree $k'+l$ , and it is a well known fact (see, for example, \cite {Ch}) that it can have at most $2(k'+l)$ zeros. Thus, our curve is $(k,l)$-regular, just as we claimed.
\smallskip

\noindent {3. {\it The plane.}}\
\smallskip

First, we observe that this map is $(k,l)$-regular in the complex sense, that is, that for any given distinct points $z_1,\ldots,z_k,w_1,\ldots w_l\in \C$, the points $\gamma(z_1),\dots,\gamma(z_k),\gamma(w_1),\ldots,\gamma(w_l)$  and complex tangent lines at last $l$ points are affinely independent over $\C$. In other words, there do not exist complex numbers $a_1,\ldots,a_k,b_1,\ldots,b_l,\xi_1,\ldots,\xi_l$, not all equal to zero, so that  the following realtions hold:
\begin {displaymath}
\sum_{i=1}^k a_i \gamma(z_i) + \sum_{j=1}^l b_j \gamma(w_j)+ \sum_{j=1}^l\xi_j \gamma'(w_j)=0,~
\sum_{i=1}^k a_i + \sum_{j=1}^l b_j=0
\end{displaymath}
The exact same proof as in the case of real open curves goes through  
when we replace real coordinates and coefficients with complex ones.\

Now, we show $(k,l)$-regularity in the real sense. Suppose therefore, towards a contradiction, that there are distinct points \mbox{$z_1,\ldots,z_k,w_1,\ldots,w_l\in \R^2$} and non-zero vectors $\xi_j \in T_{w_j}\R^2,j=1,\ldots,l$ so that $(\gamma(z_1),1),\ldots,(\gamma (z_k),1)$, $(\gamma(w_1),1),\ldots,(\gamma(w_l),1)$, $(d\gamma(\xi_1),0),\ldots,(d\gamma(\xi_l),0)$ are linearly dependent. Thus, for some real numbers $a_1,\ldots,a_k,b_1,\ldots,b_l$, not all equal to zero, we have:
\begin{displaymath}
\sum_{i=1}^k a_i \gamma(z_i) + \sum_{j=1}^l b_j \gamma(w_j)+ \sum_{j=1}^l d\gamma_{w_j}(\xi_j)=0,~
\sum_{i=1}^k a_i + \sum_{j=1}^l b_j=0
\end{displaymath}
But $d\gamma_{w}(\xi)=\xi \gamma'(w)$, where on the right side of the equation we view $\xi$ and $w$ as complex numbers and $\gamma$ as a map from $\C$ to $\C^{k+2l-1}$, and this contradicts $(k,l)$-regularity in the complex sense, which we have already established.
\proofend

\noindent {\bf{Proof of Theorem \ref{curves}}}\
We will prove that $N_{k,l}(S^1)=k+2l$ for even k. The other two claims follow immediately from (\ref {count1}) and Proposition \ref{examples}.\ 

\begin{comment}
We need to show that for $k=2k'$, $N_{k,l}(S^1) \geq 2(k'+l)$.
Suppose not. Then, $\gamma$ is a $(k,l)$-regular curve in $\R^{2(k'+l)-1}$.
Arguing as in the proof of \ref{closed}, we consider $conv(\gamma)$ and find a support hyperplane
$A^{2(k'+l)-2}$ that is tangent to $\gamma$ at $k'+l$ distinct points. What we
 want is to find an affine hyperplane $A^{2(k'+l)-2}$ that touches $\gamma$ at $l$ points and intersects it at $k=2k'$ points.

 Let $L=L^{2l-1}$ be the affine span of the tangent lines to $\gamma$ at $l$ of these points. Take $k-1$ generic points on $\gamma$ and consider  the affine
 span of $L$ and these $k-1$ points. This affine space $A$ has dimension
 $k+2l-2$, and we would be done if we can show that $A$ intersects $\gamma$ at one more point.\

 Note that the dimensions of $\gamma$ and $A$ are complimentary. Also the homological
 intersection of $\gamma$ with $A$ is zero, so if all intersection points are
 transversal then their algebraic number is zero and, in particular, the
 total number is even. In this situation $A$ touches $\gamma$ at $l$ points, and intersects $\gamma$  transversally at $2k'-1$ points. Therefore there must be another intersection point.
\proofend

2. Lower bound  $N_{k,l} (S1) \geq k+2l$ for  $k$ even.
\end{comment}

\begin{lemma} \label{curvature} Let $\gamma \subset R^n$ be a smooth curve with non-vanishing curvature vector $v$ at point $x$, and $H$ be a hyperplane, tangent to $\gamma$ at $x$ and transverse to $v$. Then, near $x$, the curve $\gamma$ doesn't cross $H$.
\end{lemma}

Proof. Let $\gamma(t)$ be parametrized by arc length  with $\gamma(0)=x$. Let \mbox{$u=\gamma'(0)$}. Then
$$\gamma(t)=x + tu + \frac{t^2}{2} v +...$$
We can assume $H$ to be the zero level hyperplane of a linear function $l$. Then, the nonvanishing curvature implies $l(v)\neq 0$, say $l(v)>0$. One has:
$$l(\gamma(t))=l(x) + tl(u) + \frac{t^2}{2}l(v) + O(t^3)=  \frac{t^2}{2}l(v) + O(t^3) > 0$$
for $t$ small enough, which proves the lemma.
\proofend
\smallskip

Now we prove the lower bound by contradiction. Assume that $S^1\subset\R^{k+2l-1}$ as a $(k,l)$-regular submanifold. In order to arrive at a contradiction, we want to find a hyperplane in $\R^{k+2l-1}$ that intersects $\gamma$ at $k$ points and touches at $l$ points, all of them distinct. Choose $l$ points $y_j$ where the curvature doesn't vanish. They must exist unless the curve $\gamma$ is straight. Let $W$ be the affine span of the tangent spaces (to $\gamma$) at these points. Let $v_j$ be the curvature vectors at $y_j$, and let $V_j$ be the affine span of $v_j$ and $W$. Note that dim $V_j \leq 2l$.

Assume that $k \geq 2$ (case $k=0$ holds because of Theorem $\ref {closed}$). Note that $\gamma$ doesn't lie in $\cup_{j=1}^l V_j$, because then there would exist an open piece of this curve contained in one of the $V_j$'s, which is impossible since $N_{k,l}(\R) \geq k+2l-1$.  Choose generic $k-1$ points $x_i \in \gamma$ that belong to neither of $V_j$ and let $U$ be the affine span of $W$ and these points. Then $U$ is a hyperplane that is transverse to $\gamma$ at points $x_i$ (by their general position) and, by Lemma $\ref{curvature}$, $\gamma$ doesn't cross $U$ at points $y_i$. Thus we have an odd number of crossings of $\gamma$ and $U$, and since $\gamma$ is closed, there must be another one, say $x_k \in U \cap \gamma$. Therefore $\gamma$ is not $(k,l)$-regular. 
\proofend

{\it E-mail address: }gordana@math.brown.edu; stojanov@math.psu.edu

\bigskip

\begin{thebibliography}{99}

\vskip 5 mm

\bibitem{A} V. Arnold. On the number of flattening points on space curves. {\it Amer. Math. Soc. Transl.} (2) {\bf 171} (1996), 11--22. 

\bibitem {Ba} A. Barvinok. {\it A course in convexity}. AMS, Providence, 2002. 

\bibitem {Bo} K. Borsuk. On the $k$-independent subsets of the Euclidean space and of the Hilbert Space. {\it Bull. Acad. Pol. Sci. Cl.III} {\bf5}, (1957),351-356.

\bibitem {BRS} V. G. Boltyansky, S. S. Ryzhkov and Yu. A. Shashkin. On $k$-regular embeddings and their applications to the theory of approximation of functions.{\it Uspekhi Mat. Nauk} {\bf15} (1960), no. 6 (96), 125--132; {\it Amer. Math. Soc. Transl.} (2) {\bf28} (1963), 211--219. 

\bibitem {Ch} E. W. Cheney. {\it Introduction to Approximation Theory}, $2^{nd}$ Ed., Chelsea Publishing Company, New York, 1996.

\bibitem {C-H} F. R. Cohen and D. Handel. $k$-regular embeddings of the plane, {\it Proc. Amer. Math. Soc.} {\bf72} (1978), 201--204. 
 
\bibitem {Gh1} M. Ghomi. Tangent bundle embeddings of manifolds in Euclidean space. Preprint.

\bibitem {Gh2} M. Ghomi. Nonexistence of skew loops on ellipsoids. {\it Proc. Amer. Math. Soc.}, in print.

\bibitem {G-S} M. Ghomi, B. Solomon. Skew loops and quadratic surfaces. {\it Comment. Math. Helv.} {\bf77} (2002), 767--782.

\bibitem {G-T} M. Ghomi, S. Tabachnikov. Totally skew embeddings of manifolds. Preprint, ArXiv: math.DG/0302288.

\bibitem {H1} D. Handel. Obstructions to 3-regular embeddings. {\it Houston J. Math} {\bf5} (1979), 339--343. 

\bibitem {H2} D. Handel. Approximation theory in the space of sections of a vector bundle. {\it Trans. Amer. Math. Soc.} {\bf 256} (1979), 383--394. 

\bibitem {H3} D. Handel. Some existence and non-existence theorems for $k$-regular maps. {\it Fund. Math.} {\bf109} (1980), 229--233.

\bibitem{H4} D. Handel. $2k$-regular maps on smooth manifolds. {\it Proc. Amer. Math. Soc.} {\bf124} (1996), 1609--1613.

\bibitem {H-S} D. Handel, J. Segal. On $k$-regular embeddings of spaces in Euclidean space. {\it Fund. Math.} {\bf106} (1980), 231--237.

\bibitem{S-S} J.-P. Sha, B. Solomon. No skew branes on non-degenerate hyperquadrics. Preprint, ArXiv: math.DG/0412197.

\bibitem{S-T} G. Stojanovic, S. Tabachnikov. Non-existence of $n$-dimensional $T$-embedded discs in ${\R}^{2n}$. Preprint, ArXiv: math.DG/0501323

\bibitem {Ta} S. Tabachnikov. On skew loops, skew branes and quadratic
hypersurfaces. {\it Moscow Math. J.} {\bf3} (2003), 681--690.

\bibitem{T-T} S. Tabachnikov, Y. Tyurina. Existence and non-existence of skew branes. Preprint, ArXiv: math.DG/0504484. 

\bibitem{V}  V. A. Vassiliev. {\it Complements of Discriminants of Smooth Maps: Topology and Applications}, Revised Edition, AMS, Providence, 1992.

\end{thebibliography}
\end{document}